\documentclass{article}
\usepackage[utf8]{inputenc}

\pdfoutput=1
\usepackage{mathtools}
\usepackage{amsfonts}
\usepackage{amssymb}
\usepackage{centernot}
\usepackage{graphicx}
\usepackage{amsthm}
\usepackage{enumerate}

\usepackage{tikz}
\usepackage{tikz-cd}
\usetikzlibrary{decorations.markings}
\usetikzlibrary{arrows}
\usetikzlibrary{calc}
\providecommand{\U}[1]{\protect\rule{.1in}{.1in}}
\newtheorem{theorem}{Theorem}

\newtheorem{corollary}[theorem]{Corollary}

\newtheorem{definition}[theorem]{Definition}

\newtheorem{lemma}[theorem]{Lemma}

\theoremstyle{remark}

\newcommand{\MQ}{\textup{MQ}}

\newcommand{\Mr}{M_A^{\textup{red}}}

\newcommand{\fr}{\phi_L^{\textup{red}}}
\newcommand{\ann}{\textup{ann}}

\allowdisplaybreaks

\begin{document}

\title{Multivariate Alexander quandles, VI. Metabelian groups and 2-component links}
\author{Lorenzo Traldi\\Lafayette College\\Easton, PA 18042, USA\\traldil@lafayette.edu
}
\date{ }
\maketitle

\begin{abstract}
We prove two properties of the modules and quandles discussed in this series. First, the fundamental multivariate Alexander quandle $Q_A(L)$ is isomorphic to the natural image of the fundamental quandle in the metabelian quotient $G(L)/G(L)''$ of the link group. Second, the medial quandle of a classical 2-component link $L$ is determined by the reduced Alexander invariant of $L$.


\end{abstract}

\section{Introduction}
We begin with a review of notation and terminology.

A (virtual) \emph{link diagram} $D$ is constructed as follows. We begin with a finite number $\mu$ of oriented circles, generically immersed in the plane. By ``generically'' we mean that the only (self-)intersections among the immersed circles are transverse double points, called \emph{crossings}. At a \emph{classical} crossing, two short subsegments are removed from one of the incident circle segments; this one is the underpassing segment, and the other is the overpassing segment. At a \emph{virtual} crossing, a small circle is drawn around the crossing. As Kauffman said \cite{vkt}, a virtual crossing is ``not really there'': we think of the two segments passing through the virtual crossing without meeting each other. The subsegments removed from the underpassing segments at classical crossings cut the original immersed circles into the \emph{arcs} of $D$. The set of arcs of $D$ is denoted $A(D)$, and the set of classical crossings of $D$ is denoted $C(D)$. A (virtual) \emph{link} $L=K_1 \cup \dots \cup K_\mu$ is an equivalence class of diagrams under the relation defined using detour moves and Reidemeister moves; each \emph{component} $K_i$ is represented by one of the $\mu$ immersed circles from which a diagram is constructed. See \cite{vkt, MI} for thorough discussions of links and their diagrams. 

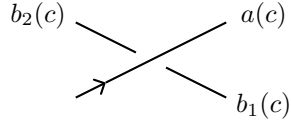
\begin{figure} [bth]
\centering
\begin{tikzpicture} [>=angle 90]
\draw [thick] (1,0.5) -- (-0.6,-0.3);
\draw [thick] [<-] (-0.6,-0.3) -- (-1,-0.5);
\draw [thick] (-1,0.5) -- (-.2,0.1);
\draw [thick] (0.2,-0.1) -- (1,-0.5);
\node at (1.5,0.6) {$a(c)$};
\node at (-1.5,0.6) {$b_2(c)$};
\node at (1.5,-0.6) {$b_1(c)$};
\end{tikzpicture}
\caption{The arcs incident at a classical crossing $c$ are indexed so that $b_1(c)$ is on the right of the overpassing arc $a(c)$.}
\label{crossfig}
\end{figure}

In the previous papers of this series we have considered several kinds of algebraic invariants of a link $L$, defined from a link diagram $D$ using generators corresponding to arcs and relations corresponding to crossings. Standard arguments verify that the algebraic structures are all invariant under detour moves and Reidemeister moves, up to isomorphism. 

\begin{definition} The \emph{fundamental quandle} $Q(L)$ is the quandle with a generator $q_a$ for each $a \in A(D)$, and a relation $q_{b_1(c)} \triangleright q_{a(c)}=q_{b_2(c)}$ for each $c \in C(D)$.
    
\end{definition} 

\begin{definition}
    
The \emph{medial quandle} $\MQ(L)$ is the medial quandle with a generator $q_a$ for each $a \in A(D)$, and a relation $q_{b_1(c)} \triangleright q_{a(c)}=q_{b_2(c)}$ for each $c \in C(D)$.
\end{definition} 

\begin{definition} \label{group}
The \emph{group} $G(L)$ is the group with a generator $g_a$ for each $a \in A(D)$, and a relation $g_{a(c)}g_{b_1(c)} g_{a(c)}^{-1}g_{b_2(c)}^{-1}=1$ for each $c \in C(D)$.
    
\end{definition}

There are natural quandle maps $Q(L) \to \MQ(L)$ and $Q(L) \to \textup{Conj} G(L)$, where $\textup{Conj} G(L)$ is the quandle defined by conjugation in $G(L)$. Here ``natural'' means that the maps are defined in the obvious ways, preserving the associations between generators and arcs. The natural map $Q(L) \to \MQ(L)$ is surjective. We refer to the image of the natural map $Q(L) \to \textup{Conj} G(L)$ as the \emph{natural image of} $Q(L)$ \emph{in} $G(L)$. Then $Q(L)$ also has a natural image in any quotient group of $G(L)$. In particular, there is a natural image of $Q(L)$ in the metabelian quotient $G(L)/G(L)''$.

Let $\Lambda_{\mu}= \mathbb Z [t_1^{\pm 1}, \dots, t_{\mu}^{\pm 1}]$ be the ring of Laurent polynomials in $\mu$ variables, with integer coefficients. Let  $\epsilon:\Lambda_\mu \to \mathbb Z$ be the \emph{augmentation map}, given by $\epsilon(t_i)=1 \thickspace \allowbreak \forall i \in \{1, \dots, \mu \} $. Then the kernel of $\epsilon$ is the \emph{augmentation ideal} of $\Lambda_{\mu}$, i.e., the ideal $I_ \mu$ generated by $\{ t_1-1, \dots, t_ \mu -1\}$. 

Given a diagram $D$ of $L=K_1 \cup \dots \cup K_\mu$, let $\kappa_D:A(D) \to \{1, \dots, \mu \}$ be the function with $\kappa_D(a)=i$ whenever $a$ is an arc of $D$ that belongs to $K_i$. Let $\Lambda_\mu^{A(D)}$ and $\Lambda_\mu^{C(D)}$ be the free $\Lambda_\mu$-modules on the sets $A(D)$ and $C(D)$, and let $\rho_D:\Lambda_\mu^{C(D)} \to \Lambda_\mu^{A(D)}$ be the $\Lambda_\mu$-linear map given by 
\[
\rho_D(c)=(1-t_{\kappa_D(b_1(c))})a(c)+t_{\kappa_D(a(c))}b_1(c) - b_2(c) \thickspace \forall c \in C(D).
\]

\begin{definition} \label{Amodule} The (multivariate) \emph{Alexander module} $M_A(L)$ is the cokernel of $\rho_D$. The canonical surjection $\Lambda_\mu^{A(D)} \to M_A(L)$ is denoted $\gamma_D$.
\end{definition}

\begin{definition} The \emph{Crowell map} is the $\Lambda_\mu$-linear map $\phi_L:M_A(L) \to I_\mu$ given by $\phi_L\gamma_D(a) = t_{\kappa_D(a)}-1 \thickspace \allowbreak \forall a \in A(D)$. The kernel of $\phi_L$ is the \emph{Alexander invariant} of $L$.
\end{definition}

We say two links $L$ and $L'$ are \emph{Crowell equivalent} if there is an isomorphism $f:M_A(L) \to M_A(L')$ with $\phi_L = \phi_{L'} f$. Notice that the definition of $\phi_L$ involves the component indices in $L$, so it is possible for $L$ and $L'$ to be inequivalent even if $L'$ is obtained from $L$ simply by re-indexing its components.

In the first paper of this series \cite{mvaq1} we introduced two quandles contained in $M_A(L)$.

\begin{definition} \label{totalalexander} Let $U(L)=\{x \in M_A(L) \mid \phi_L(x)+1 \text{ is a unit of } \Lambda_\mu\}$. Then $U(L)$ is a quandle under the operation
\[
x \triangleright y= (\phi_L(y)+1)x - \phi_L(x)y.
\]
We call this quandle the \emph{total multivariate Alexander quandle} of $L$. 

\end{definition}

\begin{definition}
    
The subquandle of $U(L)$ generated by $\gamma_D(A(D))$ is the \emph{fundamental multivariate Alexander quandle} of $L$, denoted $Q_A(L)$.

\end{definition}

We are now ready to state the first theorems of the present paper.

\begin{theorem}
\label{metabelian}
For any link $L$, $U(L)$ is isomorphic to $\textup{Conj}(G(L)/G(L)'')$, the conjugation quandle of the metabelian quotient of $G(L)$. Moreover, there is an isomorphism $U(L) \cong \textup{Conj}(G(L)/G(L)'')$ that maps $Q_A(L)$ isomorphically onto the natural image of $Q(L)$ in $G(L)/G(L)''$.
\end{theorem}

Theorem \ref{metabelian} leads to the following refinement of \cite[Theorem 11]{mvaq1}.

\begin{theorem}
\label{main}
Consider the following statements about links $L=K_1 \cup \dots \cup K_{\mu}$ and $L'=K'_1 \cup \dots \cup K'_{\mu'}$. 
\begin{enumerate}
\item $Q(L) \cong Q(L')$.
\item $Q_A(L) \cong Q_A(L')$.
\item The components of $L$ and $L'$ can be re-indexed so that the resulting links are Crowell equivalent.
\item $G(L)/G(L)'' \cong G(L')/G(L')''$.
\item $U(L) \cong U(L')$.
\end{enumerate}
The implications $1 \implies 2 \implies 3 \implies 4 \implies 5$ all hold for virtual links, and the converses all fail even when restricted to classical links.
\end{theorem}

Let $\Lambda$ be the ring $\Lambda= \mathbb Z [t^{\pm 1}]$ of Laurent polynomials in one variable, with integer coefficients, and let $\tau:\Lambda_ \mu \to \Lambda$ be the ring homomorphism with $\tau(t_i)=t$ $\forall i$. Then the reduced (one-variable) Alexander module $M_A^{\text{red}}(L)$ is the $\Lambda$-module obtained by applying $\tau$ to the definition of $M_A(L)$ given above. Equivalently, $\Mr(L)$ is the tensor product $M_A(L) \otimes_{\Lambda_\mu}\Lambda$, where $\Lambda$ is considered a $\Lambda_\mu$-module via $\tau$. For knots, of course, $\Mr(L)=
M_A(L)$.

We use $\varsigma_D(a)$ to denote the generator of $\Mr(L)$ corresponding to an arc $a \in A(D)$. When we think of $\Mr(L)$ as $M_A(L) \otimes_{\Lambda_\mu}\Lambda$, $\varsigma_D(a)$ is $\gamma_D(a) \otimes 1$.

The tensor product of the Crowell map with the identity map of $\Lambda$ is a map $\phi_\tau: \Mr(L) \to I_\mu \otimes _{\Lambda_\mu}\Lambda$. As discussed in \cite{mvaq4}, it is not hard to describe $\phi_\tau$ explicitly. The tensor product $I_\mu \otimes _{\Lambda_\mu}\Lambda$ is isomorphic to the direct sum $\Lambda \oplus \mathbb Z ^{\mu-1}$. Here $\mathbb Z ^{\mu-1}$ is the free abelian group of rank $\mu-1$, considered as a trivial $\Lambda$-module. (That is, the scalar multiplication is given by $tx=x \thickspace \allowbreak \forall x \in \mathbb Z ^ {\mu-1}$.) With this isomorphism in hand, $\phi_\tau$ may be described as follows: If $a \in A(D)$ and $\kappa_D(a)=1$, then $\phi_\tau(\varsigma_D(a)) = (1,0,\dots,0) \in \Lambda \oplus \mathbb Z ^{\mu-1}$; and if $a \in A(D)$ and $\kappa_D(a)=i>1$, then $\phi_\tau(\varsigma_D(a)) = (1,0,\dots,0,1,0,\dots,0) \in \Lambda \oplus \mathbb Z ^{\mu-1}$, with the second $1$ in the $i$th coordinate. 

The first coordinate of $\phi_\tau:\Mr(L) \to \Lambda  \oplus \mathbb Z ^{\mu-1}$ is an epimorphism $\Mr(L) \to \Lambda$, under which $\varsigma_D(a) \mapsto 1 \thickspace \allowbreak \forall a \in A(D)$. This is the \emph{reduced Crowell map}, denoted $\fr$. The kernel of $\fr$ is the \emph{reduced Alexander invariant} of $L$. An epimorphism onto a free module must split, so $\Mr(L) \cong \ker \fr \oplus \Lambda$. Hence $\Mr(L)$ is determined up to isomorphism by $\ker \fr$.
    
Quandles were introduced into knot theory by Joyce \cite{J} and Matveev \cite{M}. Both authors discussed the fact that the (reduced) Alexander invariant and module of a classical knot yield quandles under the operation $x \triangleright y = tx + (1-t)y$. Indeed, Joyce's Theorems 17.2 and 17.3 have the following corollaries. (Joyce called $\MQ(L)$ the ``abelian quandle of $L$'' and denoted it $\textup{AbQ}(L)$.)

\begin{corollary} \label{det1} \cite{J} If two classical knots have isomorphic medial quandles, then they also have isomorphic reduced Alexander invariants.
\end{corollary}
\begin{corollary} \label{det2} \cite{J} If two classical knots have isomorphic reduced Alexander invariants, then they also have isomorphic medial quandles.
\end{corollary}

Notice that Corollaries \ref{det1} and \ref{det2} are converses of each other. We state them separately because they extend to virtual links in different ways \cite{mvaq4, mvaq5}: Corollary \ref{det1} holds for all virtual links, but examples show that Corollary \ref{det2} fails for virtual links with $\mu \geq 2$, and also fails for classical links with $\mu \geq 3$. Here we complete the picture by showing that in contrast to the virtual case, Corollary \ref{det2} holds for classical links with $\mu = 2$. As Corollary \ref{det1} also holds for these links, we deduce the following.  

\begin{theorem}
\label{maintwo}
Two classical 2-component links have isomorphic medial quandles if and only if they have isomorphic reduced Alexander invariants.
\end{theorem}

\section{Theorem \ref{metabelian}}
\label{section2}
Theorem \ref{metabelian} follows from two algebraic results stated by Crowell \cite{C1, C3}. To set up these results, let $f:G\to H$ be a group epimorphism with kernel $K$, and let $K'$ be the commutator subgroup of $K$. Then $K/K'$ is an $H$-module with the scalar multiplication defined using conjugation in $G$. That is, if $h=f(g) \in H$ and $k \in K$ then $h \cdot kK' = gkg^{-1}K'$. Let $\mathbb Z f:\mathbb Z G \to \mathbb Z H$ be the linear extension of $f$. Let $IG$ be the augmentation ideal of $G$, i.e., the kernel of the augmentation map $\epsilon: \mathbb Z G \to \mathbb Z$ given by $\epsilon (\sum n_i g_i) = \sum n_i$. Define $IH$ in the same way.

\begin{theorem} \label{seq} \cite{C1,C3}
In this situation there is an exact sequence of $H$-modules
\begin{equation*}
0 \to K/K' \xrightarrow{\psi} \mathbb Z H \otimes_{G} IG \xrightarrow{\phi} IH \to 0 \text{,}
\end{equation*}
where $\mathbb Z H \otimes_{G} IG$ is an $H$-module via multiplication in $\mathbb Z H$, $\psi$ is given by $\psi (kK')=1 \otimes (k-1)$, and $\phi$ is given by $\phi(x \otimes y)=x \cdot \mathbb Z f(y)$.
\end{theorem}

Crowell \cite{C1, C3} discussed four different descriptions of the $H$-module given by the tensor product $\mathbb Z H \otimes_{G} IG$. Here they are.

\begin{enumerate}
\item $\mathbb Z H \otimes_{G} IG$ is characterized up to isomorphism by a universal property involving crossed homomorphisms. 
\item $\mathbb Z H \otimes_{G} IG$ is isomorphic to the quotient module $IG/(IG \cdot \ker \mathbb Z f)$. 
\item Let $i:G \to X$ be a one-to-one correspondence between $G$ and a set $X$. Then $\mathbb Z H \otimes_{G} IG$ is isomorphic to the quotient of the free $H$-module $(\mathbb Z H)^X$ by the submodule generated by the elements $i(g_1g_2)-i(g_1)-f(g_1)i(g_2) $ with $g_1,g_2 \in G$. 
\item Let $<S;R>$ be a presentation of $G$. Then $\mathbb Z H \otimes_{G} IG$ is isomorphic to the quotient of the free $H$-module $(\mathbb Z H)^S$ by the submodule $N$ generated by elements represented by the rows of $\mathbb Z f(J)$, where $J$ is a Jacobian matrix derived from the presentation $<S;R>$ using the free differential calculus of Fox.
\end{enumerate}

We do not reference Crowell \cite{C1, C3} for Theorem \ref{seq} to give the impression that he discovered it. In fact, Theorem \ref{seq} is a standard result, and we do not know who stated it first. Some version of the exact sequence appears in every text that presents the theory of group cohomology. For instance, Cartan and Eilenberg \cite[Sec.\ XIV.4]{CE} mentioned the sequence with the second different description of $\mathbb Z H \otimes_{G} IG$, Hilton and Stammbach \cite[Sec.\ VI.6]{HS} mentioned the sequence with the module $\mathbb Z H \otimes_{G} IG$, and Mac Lane \cite[Sec.\ IV.6]{Mac} mentioned the sequence with the third different description of $\mathbb Z H \otimes_{G} IG$. As Theorem \ref{seq} is a standard result, we do not present a proof.

One reason we reference Crowell \cite{C1,C3} for Theorem \ref{seq} is that he provided a detailed discussion of the fact that the $H$-module $\mathbb Z H \otimes_{G} IG$ has the four different descriptions listed above. In particular, there is an isomorphism between $\mathbb Z H \otimes_{G} \nobreak IG$ and the module $(\mathbb Z H)^S/N$ of the fourth different description defined by $1 \otimes (g_s-1) \mapsto s+N$, where $g_s$ is the element of $G$ represented by an element $s \in S$. This isomorphism will be useful in our proof of Theorem \ref{metabelian}. 

Another reason we reference Crowell for Theorem \ref{seq} is that he proved the next result, which we have not seen stated in any standard discussion of group cohomology. 

\begin{theorem} \label{Cgroup} \cite{C1}
Consider the exact sequence of Theorem \ref{seq}, and let $
U(\phi)=\{x \in \mathbb Z H \otimes_{G} IG \mid 1+\phi(x) \text{ is a unit in }\mathbb Z H \}$.
\begin{enumerate}
\item \label{parta} $U(\phi)$ is a group under the operation $x \cdot y = x + (1+\phi(x))y$.
\item \label{partb} There is an isomorphism $\eta:G/K' \to U(\phi)$ given by $\eta(gK')= 1 \otimes (g-1)$.
\end{enumerate}
\end{theorem}

\begin{proof} For the reader's convenience we sketch a proof adapted from \cite {C1}. For part \ref{parta}, the reader can easily check that the operation $\cdot$ is associative, that $0$ plays the role of the identity, and that the inverse of an element $x$ is $-(1+\phi(x))^{-1}x$. 

For part \ref{partb}, note first that $\eta$ is well defined: If $g_1K'=g_2K'$ then
\begin{align*}
0 &=\psi(g_1g_2^{-1}K')= 1 \otimes(g_1g_2^{-1}-1)=1 \otimes(g_1g_2^{-1}-g_1)+1 \otimes(g_1-1)
\\
&=1 \otimes(g_1g_2^{-1}(1-g_2))+1 \otimes(g_1-1)=f(g_1g_2^{-1}) \otimes(1-g_2)+1 \otimes(g_1-1)
\\
&=1 \otimes (1-g_2)+1 \otimes(g_1-1) = -1 \otimes (g_2-1)+1 \otimes(g_1-1) \text{,}
\end{align*}
so $1 \otimes(g_1-1)=1 \otimes(g_2-1)$.

Also, $\eta$ is a homomorphism: If $g_1,g_2 \in G$ then 
\begin{align*}
\eta(g_1K' \cdot g_2K')&=\eta(g_1g_2K')=1 \otimes(g_1g_2-1) = 1 \otimes(g_1-1)+1 \otimes(g_1g_2-g_1)
\\
 &= 1 \otimes(g_1-1)+1 \otimes(g_1(g_2-1))
 \\
 &=1 \otimes(g_1-1)+f(g_1) \otimes(g_2-1)
\\
&=1 \otimes(g_1-1)+(1+f(g_1)-1)(1 \otimes (g_2-1))
\\
&=1 \otimes(g_1-1)+(1+\phi(1 \otimes (g_1-1)))(1 \otimes (g_2-1))
\\
&=(1 \otimes(g_1-1)) \cdot (1 \otimes (g_2-1))=\eta(g_1K') \cdot \eta (g_2K').
\end{align*}

To show that $\eta$ is surjective, recall that the units of the group ring $\mathbb Z H$ are the elements of the form $\pm h, h \in H$. If $x \in \mathbb Z H \otimes_{G} IG$ and $1+\phi(x)$ is a unit in $\mathbb Z H$, then $\epsilon(1+\phi(x))=1+\epsilon(\phi(x))=1+0=1$, so $1+\phi(x)$ must be a unit of the form $+h$, not $-h$. Let $1+\phi(x)=h=f(g)$. Then $\phi(x) = f(g)-1=\phi(1 \otimes (g-1))$, so $x-(1 \otimes (g-1)) \in \ker \phi$. According to Theorem \ref{seq}, it follows that there is a $k \in K$ with $x-(1 \otimes (g-1)) = \psi(kK')=1 \otimes (k-1)$. Then 
\begin{align*}
\eta ((kg)K')&=1 \otimes (kg-1) = f(k^{-1}) \otimes (kg-1) = 1 \otimes (k^{-1}(kg-1))
\\
&  = 1 \otimes (g-k^{-1})= 1 \otimes (g-1)+ 1 \otimes (1 - k^{-1}) 
\\
&= 1 \otimes (g-1)+ 1 \otimes ( k^{-1}(k-1))
\\
&= 1 \otimes (g-1)+ f(k^{-1}) \otimes (k-1) = 1 \otimes (g-1)+ 1 \otimes (k-1) = x.
\end{align*}

It remains to show that $\eta$ is injective. If $gK' \in \ker \eta$ then $1 \otimes (g-1)=0$, so certainly $0=\phi(1 \otimes (g-1)) = f(g)-1 \in IH$. Then $f(g)=1 \in H$, so $g \in \ker f = K$. As $\psi(gK') = 1 \otimes (g-1) =0$, it follows that $gK' \in \ker \psi$. But Theorem \ref{seq} tells us that $\psi$ is injective, so it must be that $gK' = 1 \in K/K'$. Then $g \in K'$, so $gK' =1$ also in $G/K'$.
\end{proof}

\begin{corollary} \label{Cgroupcor}
In the situation of Theorems \ref{seq} and \ref{Cgroup}, there is a quandle structure on $U(\phi)$ given by
\[
x \triangleright y = y+(1+\phi(y))x-(1+\phi(y))(1+\phi(x))(1+\phi(y))^{-1}y.
\]
An isomorphism between this quandle and the conjugation quandle $\textup{Conj}(G/K')$ is defined by $\eta$. 
\end{corollary}
\begin{proof}
The conjugation quandle of the group $U(\phi)$ specified in part \ref{parta} of Theorem \ref{Cgroup} is defined by the operation
\begin{align*}
x \triangleright y &= y \cdot xy^{-1} = y + (1+\phi(y))xy^{-1}
=y+(1+\phi(y))(x+(1+\phi(x))y^{-1})
\\
&=y+(1+\phi(y))x+(1+\phi(y))(1+\phi(x))(-(1+\phi(y))^{-1})y.
\end{align*}
Part \ref{partb} of Theorem \ref{Cgroup} tells us that $\eta$ maps this quandle isomorphically onto $\textup{Conj}(G/K')$. 
 \end{proof}

We are now ready to prove Theorem \ref{metabelian}. Given a link $L$ with a diagram $D$, we apply the machinery discussed in this section to the canonical epimorphism $f:G(L) \to H$, where $H$ is the abelianization $G(L)/G(L)'$. Then $K=\ker f = G(L)'$, $K'=G(L)''$, and $\mathbb Z H$ may be identified with $\Lambda _ \mu$, with $f(g_a)$ corresponding to $t_{\kappa_D(a)}$ for each $a \in A(D)$. The fourth different description of $\mathbb Z H \otimes_{G} IG$ is the $H$-module $\mathbb Z H ^{A(D)}/N$, where $N$ is the submodule of $\mathbb Z H ^{A(D)}$ generated by elements represented by the rows of a Jacobian matrix $\mathbb Z f(J)$ obtained by applying Fox's free differential calculus to Definition \ref{group}. This description of $\mathbb Z H ^{A(D)}/N$ is the same as the definition of $M_A(L)$ in Definition \ref{Amodule}, because each crossing $c \in C(D)$ has $\rho_D(c)$ equal to the element of $\mathbb Z H ^{A(D)}$ represented by the row of $\mathbb Z f (J)$ corresponding to the relator $g_{a(c)}g_{b_1(c)} g_{a(c)}^{-1}g_{b_2(c)}^{-1}$ of Definition \ref{group}. 

As discussed in the penultimate paragraph before the statement of Theorem \ref{Cgroup}, it follows that there is an isomorphism $M_A(L) \to \mathbb Z H \otimes_{G} IG$ given by $\gamma_D(a) \mapsto 1 \otimes (g_a-1) \thickspace \allowbreak \forall a \in A(D)$. Notice that the Crowell map $\phi_L$ corresponds to the map $\phi$ of Theorem \ref{seq} under this isomorphism: if $a \in A(D)$ then $\phi_L(\gamma_D(a)) = t_{\kappa_D(a)}-1=f(g_a)-1=\phi(1 \otimes(g_a-1))$. The image of the quandle $U(L)$ under this isomorphism is the quandle $U(\phi)$ mentioned in Corollary \ref{Cgroupcor}, and the image of $Q_A(L)$ is the subquandle of $U(\phi)$ generated by the elements $1 \otimes (g_a-1)$ with $a \in A(D)$. Composing with the isomorphism $\eta$ of Theorem \ref{Cgroup} and Corollary \ref{Cgroupcor}, we obtain a quandle isomorphism $\textup{Conj} 
(G(L)/G(L)'') \to U(L)$ under which $g_a G(L)'' \mapsto \gamma_D(a) \thickspace \allowbreak \forall a \in A(D)$. The latter property guarantees that the isomorphism maps the natural image of $Q(L)$ in $\textup{Conj} 
(G(L)/G(L)'')$ onto $Q_A(L)$.

\section{Theorem \ref{main}}
\label{section3}

The implications $1  \implies 2$ and $2 \implies 3$ of Theorem \ref{main} are proven just as in \cite{mvaq1}. The implication $3 \implies 4$ follows from Theorem \ref{Cgroup}, and $4 \implies 5$ follows from Theorem \ref{metabelian}.

Examples given in \cite{mvaq1} show that $2 \centernot \implies 1$ and $4 \centernot \implies 3$. The former are classical knots with Alexander polynomial 1, and the latter are classical links with homeomorphic complements and distinct Alexander polynomials. The classical links $W$ and $7^2_8$ analyzed in \cite{mvaq3} demonstrate that $3 \centernot \implies 2$. 

For $5 \centernot \implies 4$, it is well known that the unknot $K_0$ and the Hopf link $H$ have $G(K_0)/G(K_0)'' \cong G(K_0) \cong \mathbb Z$ and $G(H)/G(H)'' \cong G(H) \cong \mathbb Z \oplus \mathbb Z$. The conjugation quandles of these two groups are isomorphic: both are countably infinite and trivial. (That is, they satisfy $x \triangleright y \equiv x$.)

\section{Theorem \ref{maintwo}}
\label{section4}
 
Proving Theorem \ref{maintwo} requires ideas and results from \cite{periadd, peri, mvaq5}. We review these ideas and results very briefly, and refer to the original papers for detailed discussions. 

As mentioned in the introduction, the reduced Alexander module $\Mr(L)$ is defined by applying $\tau$ to the definition of $M_A(L)$ in Definition \ref{Amodule}; the module generator corresponding to an arc $a \in A(D)$ is denoted $\varsigma_D(a)$. The maps $\phi_\tau:\Mr(L) \to \Lambda \oplus \mathbb Z ^{\mu-1}$ and $\fr:\Mr(L) \to \Lambda$ are also defined in the introduction.

A peripheral structure for the reduced Alexander module $\Mr(L)$ was introduced in \cite{peri}. If $L=K_1 \cup \dots \cup K_ \mu$ then each component $K_i$ has a corresponding set of meridians $M_i(L) \subset \Mr(L)$, and a single corresponding longitude $\chi_i(L) \in \Mr(L)$. Here are six important properties, numbered for later reference. The first five are from \cite{peri} and the sixth is from \cite{periadd}.
\begin{enumerate}
\setcounter{enumi}{0}
\item \label{prop1} If $D$ is a diagram of $L$ then for each $a \in A(D)$, $\varsigma_D(a) \in M_{\kappa_D(a)}(L)$. 
\item \label{prop2} If $1 \leq i \leq \mu$, then $\phi_\tau$ is constant on $M_i(L)$. 
\item \label{prop3} If $1 \leq i \leq \mu$, then $\fr(\chi_i(L)) = 0$. 
\item \label{prop4} The longitudes generate the submodule of $\Mr(L)$ annihilated by $t-1$.
\item \label{prop5} The kernel of $\phi_\tau$ is $(1-t) \ker \fr$.
\item \label{prop6} If $L$ is a classical link, then $\sum_{i=1}^ \mu \chi_i(L)=0$.
\end{enumerate}

In \cite{mvaq5}, the following definition was introduced.

\begin{definition}
\label{qquandle}
Suppose $M$ is a $\Lambda$-module with a submodule $N$, $I$ is a nonempty set, $m_i \in M\thickspace \allowbreak \forall i \in I$, and $m_i-m_j \in N \thickspace \allowbreak \forall i,j \in I$. Suppose also that for each $i \in I$, $X_i$ is a submodule of $N$ such that $(1-t) \cdot X_i = 0$. Then there is an associated medial quandle $Q(N,(m_i), \allowbreak (X_i))$, defined as follows.
\begin{enumerate} [(a)]
    \item For each $i \in I$, let $Q_i=N/X_i$. If $i \neq j \in I$ then the sets $Q_i$ and $Q_j$ are understood to be disjoint. Define the set $Q(N,(m_i), \allowbreak (X_i))$ by
    \[
    Q(N,(m_i), \allowbreak (X_i))=\bigcup_{i \in I} Q_i .
    \]
    \item Define an operation $\triangleright$ on $Q(N,(m_i), \allowbreak (X_i))$ as follows. If $i,j \in I$, $x\in Q_i$ and $y \in Q_j$, then
    \[
    x \triangleright y = m_j-m_i+tx+(1-t)y+X_i \in Q_i.
    \]
\end{enumerate}
\end{definition}

Here are four important properties from \cite{mvaq5}, numbered to make it easy to reference them together with the six properties from \cite{periadd, peri} mentioned above. 
\begin{enumerate} 
\setcounter{enumi}{6}
\item \label{prop7} Every medial quandle can be realized as $Q(N,(m_i), \allowbreak (X_i))$ for appropriate choices of $M,N, I, (m_i)$ and $(X_i)$.
\item \label{prop8} For any fixed element $m \in M$, the quandle $Q(N,(m_i), \allowbreak (X_i))$ is not changed if every $m_i$ is replaced with $m_i - m$.
\item \label{prop9} For any choice of elements $n_i \in N$, the quandles $Q(N,(m_i), \allowbreak (X_i))$ and $Q(N,(m_i+(1-t)n_i), \allowbreak (X_i))$ are isomorphic.
\item \label{prop10} If $L=K_1 \cup \dots \cup K_ \mu$, $I=\{1, \dots , \mu \}$, $M=\Mr(L)$, $N= \ker \fr$, each $m_i$ is an arbitrary element of $M_i(L)$, and each $X_i$ is the cyclic submodule of $\ker \fr$ generated by $\chi_i(L)$, then $\MQ(L) \cong Q(N,(m_i), \allowbreak (X_i))$.
\end{enumerate}

\begin{lemma} \label{annlem} Let $L=K_1 \cup K_2$ be a classical link of two components. Then $\chi_1(L)=-\chi_2(L)$, and the cyclic submodule of $\Mr(L)$ generated by either of the two longitudes is $X_1=X_2=\ann(1-t) = \{x \in \Mr(L) \mid (1-t)x=0 \}$.
\end{lemma}
\begin{proof}
The lemma follows from properties \ref{prop4} and \ref{prop6}.
\end{proof}

We say a $\Lambda$-module $M$ is trivial if it has $tx=x \thickspace \allowbreak \forall x \in M$.

\begin{lemma}
\label{quotient}
Let $L=K_1 \cup K_2$ be a classical link of two components. Then the quotient module $\ker \fr / (1-t) \ker \fr $ is isomorphic to the trivial $\Lambda$-module $\mathbb Z$. If $m_1 \in M_1(L)$ and $m_2 \in M_2(L)$, then $m_1-m_2 \in \ker \fr$ and the coset $m_1-m_2+(1-t)\ker\fr$ generates $\ker \fr / (1-t) \ker \fr $.
\end{lemma}
\begin{proof}
As discussed in the introduction, there is a $\Lambda$-linear epimorphism $\phi_\tau:\Mr(L) \to \Lambda \oplus \mathbb Z$, where $\mathbb Z$ is trivial as a $\Lambda$-module. The first coordinate of $\phi_\tau$ is the $\Lambda$-linear epimorphism $\fr:\Mr(L) \to \Lambda$. It follows that the second coordinate of $\phi_\tau$ restricts to a $\Lambda$-linear epimorphism $\widetilde \phi:\ker \fr \to \mathbb Z$, and of course this epimorphism has the property that $\mathbb Z \cong \ker \fr / \ker \widetilde \phi$. As $\ker \widetilde \phi = \ker \phi_\tau$, property \ref{prop5} tells us that $\ker \widetilde \phi = (1-t) \ker \fr$.

If $m_1 \in M_1(L)$ and $m_2 \in M_2(L)$, then properties \ref{prop1} and \ref{prop2}, along with the definition of $\phi_\tau$ in the introduction, tell us that $\phi_\tau(m_1-m_2)=(1,0)-(1,1)=(0,-1) \in \Lambda \oplus \mathbb Z$. Therefore $m_1-m_2 \in \ker \fr$ and $\widetilde \phi (m_1-m_2) = -1 \in \mathbb Z$, so $\widetilde \phi (m_1-m_2)$ generates $\mathbb Z$. It follows that $m_1-m_2+\ker \widetilde \phi$ generates the quotient $\ker \fr / \ker \widetilde \phi$.
\end{proof}

\begin{theorem}
\label{inv}
Let $L=K_1 \cup K_2$ be a classical link of two components, let $X=\ann(1-t)$ be the submodule of $\ker \fr$ annihilated by $1-t$, and let 
\[
D=\{ d \in \ker \fr \mid d + (1-t)\ker \fr \textup{ generates } \ker \fr / (1-t) \ker \fr  \}.
\]
Then every $d \in D$ has $Q(\ker \fr, (d,0),(X,X)) \cong \MQ(L)$.
\end{theorem}
\begin{proof} Suppose $d \in D$.

To begin, note that according to property \ref{prop8}, the quandles $Q(\ker \fr, (d,0), \allowbreak (X,X)$ and $Q(\ker \fr, (0,-d),(X,X))$ are isomorphic. Permuting the index set $I$ has no effect on the quandle $Q(N,(m_i),(X_i))$ of Definition \ref{qquandle}, so we conclude that the quandles $Q(\ker \fr, \allowbreak (d,0), (X,X))$ and $Q(\ker \fr, (-d,0), \allowbreak (X,X))$ are isomorphic. That is, we may replace $d$ with $-d$ if we like, without changing the quandle $Q(\ker \fr, \allowbreak (d,0), (X,X))$.

Let $m_1 \in M_1(L)$ and $m_2 \in M_2(L)$. Then Lemma \ref{quotient} tells us that the quotient module $\ker \fr / (1-t) \ker \fr $  is isomorphic to $\mathbb Z$, and is generated by the coset $m_1-m_2+(1-t)\ker\fr$. As $d + (1-t)\ker \fr$ also generates $\ker \fr / (1-t) \ker \fr $, it follows that 
\[
m_1-m_2+(1-t)\ker\fr = \pm d + (1-t)\ker \fr.
\]
According to the preceding paragraph, we may ignore the $\pm$ without loss of generality, by replacing $d$ with $-d$ if necessary.

Then $m_1-m_2+(1-t)\ker\fr = d + (1-t)\ker \fr$, so $m_1-m_2-d$ is an element of $ (1-t) \ker \fr$. According to property \ref{prop9}, it follows that $Q(\ker \fr, \allowbreak (d,0), (X,X))$ is isomorphic to 
\[
Q(\ker \fr, \allowbreak (d+m_1-m_2-d,0), (X,X))
=Q(\ker \fr, \allowbreak (m_1-m_2,0), (X,X)).
\]
Applying property \ref{prop8}, we conclude that $Q(\ker \fr, \allowbreak (d,0), (X,X))$ is isomorphic to 
\[
Q(\ker \fr, \allowbreak (m_1-m_2+m_2,0+m_2), (X,X))
=Q(\ker \fr, \allowbreak (m_1,m_2), (X,X)).
\]

Lemma \ref{annlem} tells us that $X=X_1=X_2$, and property \ref{prop10} tells us that $Q(\ker \fr, (m_1,m_2),(X_1,X_2)) \cong \MQ(L)$. \end{proof}

We are now ready to prove that Corollary \ref{det2} holds for classical 2-component links. Suppose $f:\ker \fr \to \ker \phi_{L'}^{\textup{red}}$ is an isomorphism, where $L$ and $L'$ are classical 2-component links. Let $X'$ and $D'$ be the subsets of $\ker \phi_{L'}^{\textup{red}}$ defined as in Theorem \ref{inv}. Then the fact that $f$ is an isomorphism implies $f(X)=X'$ and $f(D)=D'$. The fact that $f$ is an isomorphism also implies that for any $d \in D$, $Q(\ker \fr, \allowbreak (d,0), (X,X))$ is isomorphic to $Q(\ker \phi_{L'}^{\textup{red}}, \allowbreak (f(d),0), (f(X),f(X)))$.  Then Theorem \ref{inv} tells us
\begin{align*}
\MQ(L) &\cong Q(\ker \fr, \allowbreak (d,0), (X,X)) \cong Q(\ker \phi_{L'}^{\textup{red}}, \allowbreak (f(d),0), (f(X),f(X)))
\\ & = Q(\ker \phi_{L'}^{\textup{red}}, \allowbreak (f(d),0), (X',X')) \cong \MQ(L').
\end{align*}

\end{document}